\documentclass[11pt,reqno]{amsart}
\usepackage{amssymb}
\usepackage{amsmath}

 \newtheorem{thm}{Theorem}[section]
 
 \newtheorem{lem}[thm]{Lemma}

 \theoremstyle{definition}
 
 \theoremstyle{remark}
 
 \numberwithin{equation}{subsection}

\begin{document}

\title[Eigenvalues and lambda constants]
{Eigenvalues and lambda constants on Riemannian submersions}

\author{ Li Ma, Anqiang Zhu}
\address{Li Ma, Department of Mathematical Sciences, Tsinghua University,
 Peking 100084, P. R. China}

\email{lma@math.tsinghua.edu.cn}

\thanks{The research is partially supported by the National Natural Science
Foundation of China 10631020 and SRFDP 20060003002. The first
named author would like to thank Prof. S.Donaldson for reminding
him the relation between the Bochner formula and Perelman's
F-functional on spin manifolds.}

\begin{abstract}
Given a Riemannian submersion, we study the relation between
lambda constants introduced by G.Perelman on the base manifold and
the total space of a Riemannian submersion. We also discuss the
relationship between the first eigenvalues of Laplacians on the
base manifold and that of the total space. The quantities on
warped products are discussed in detail.


\end{abstract}

\date{}

\dedicatory{}

\maketitle
\section{introduction}
The aim of this paper is to study the relationship between Lambda
constants, which will be defined below, the first eigenvalues of
Laplacians on the base manifold and that of the total space of a
Riemannian submersion. Even though these are classical objects in
Differential Geometry and in Mathematical Analysis (see
\cite{com}), they are missed in literature.

Motivated by Physical background and the Gross Log-Sobolev
inequality, G. Perelman \cite{perelman} introduced the
F-functional
$$
F(g,f)=\int_{M}e^{-f}(|\nabla f|^{2}+R)dv_g,$$ where $R$ is the
scalar curvature of the Riemannian manifold $(M,g)$ ( see
\cite{besse}) and $dv_g$ is the volume element of the metric $g$,
with its infimum defined by
$$ \lambda(g)=inf\{F(g,f):f\in
C^{\infty}_{c}(M),\int_{M}e^{-f}dv_g=1\};
$$
which will be called the \emph{lambda constant}. It has been
showed in \cite{Le} and \cite{F06} that there is a closed
relationship between the Yamabe constants and the lambda
constants. Recall that the Yamabe constant $Y(g)$ on $(M,g)$ is
the infimum of the action functional
$$ A(g)=\int_M Rdv_g,
$$
which was introduced by D.Hilbert, among the conformal class $[g]$
with fixed unit volume.

 Let us recall some basic notations and concepts about Riemannian
submersions from \cite{besse}. Let $p:M\to B$ be a Riemannian
submersion with compact fiber $F$. Put $F_b=p^{-1}(b)$ for $b\in
B$. Given a smooth function $u^M$ on $M$, we define a smooth
function $u^B$ on $B$ by
$$ p_{\star}(u^Mdv_M)=u^Bdv_B.
$$
In other word, we have
$$
u^B=\int_Fudv_F.
$$
Let $N$ be the mean curvature vector field to the fibers $F$. Let
$A$ be the curvature of the horizontal distribution and let $T$ be
the second fundamental form of the fibers $F$. We denote by
$\nabla^M$ and $\nabla^B$ the covariant derivatives on $M$ and on
$B$ respectively. We always denote by $G^M$, $G^F$, and $G^B$ the
same kind of geometrical quantities on $M$, on $F$, and on $B$
respectively. For example, $\Delta^M$ the Laplacian operator on $M$.
We denote by $\nabla_{hor}$ the horizontal part of the covariant
derivative on $M$. That is
$$
\nabla_{hor}u=\nabla_{e_{\alpha}}ue_{\alpha}
$$
where $\{e_{\alpha}\}$ is an orthonormal basis of $T_{hor}M$ at a
point $m$. We define
$$
\breve{\delta}N(m)=-\sum_{\alpha}<\nabla_{e_{\alpha}}^MN,e_{\alpha}>.
$$
Then we have
$$
R^M=R^B+R^F-|A|^2-|T|^2-|N|^2-2\breve{\delta}N
$$
and
$$
\Delta^Mu=\nabla^2_{hor}u+\Delta^Fu-<\nabla_{hor}u,N>.
$$

When $M=B\times F$ with the warped metric $g=g^B+e^{2f}g^F$, where
$f\in C^{\infty}(B)$, we have
$$
\Delta^Mu=\Delta^Bu+e^{-2f}\Delta_{g^F}u+k<\nabla^Bu,\nabla^Bf>.
$$
where $k=dim(F)$.

Let
$$
A_0=\frac{1}{2}inf_M(-\sum\nabla_{e_{\alpha}}<N,e_{\alpha}>)
$$
and
$$
A_1=\frac{1}{2}sup_M(-\sum\nabla_{e_{\alpha}}<N,e_{\alpha}>).
$$

For the comparison of first eigenvalues on the base and total
space, we have the following result.
\begin{thm}\label{thm1} Given a compact Riemannian submersion $(M,g)$.
 Then we have the following inequality
$$
\lambda_1^M\leq \lambda_1^B+\bar{\lambda}_F+A_1.
$$
Here
$$
\bar{\lambda}_F=sup_{b\in B}\lambda_1^F(b).
 $$
 Furthermore, when $(M,g)$ is a warped product, we also have
 $$
\lambda_1^B+A_0\leq \lambda_1^M.
 $$
\end{thm}
Since the Euler-Lagrange operator of the F-functional $F(g,f)$ is
of the same type of the Laplace operator, we can derive the
following result
\begin{thm}\label{thm2}
Suppose $(M,g^{M})$ is a warped product of $(B,\breve{g})$ and
$(F,g_{0})$, $g^{M}=\breve{g}+f^{2}g_{0}$, where $f$ is a function
on $B$. If the scalar curvature $\hat{R_{0}}$ of the fiber
$(F,g_{0})$ is bounded from below such that $\hat{R_{0}}\geq
r_{0}$, then we have
$$
\lambda^{M}\geq\lambda^{B}+c,
$$
where $ c=inf (
2p\frac{\Delta^{B}f}{f}+p(7-p)\frac{|df|^{2}}{f^{2}}+\frac{1}{f^{2}}r_{0}).$
\end{thm}
One may ask if the similar result is true for W-functional. Recall
the quantity in question is $\mu(g,\tau)$, which is defined by
$$
\mu(g,\tau)=inf \{W(g,f,\tau)|f\in
C^{\infty}_c,\tau>0,\frac{1}{(4\pi\tau)^{\frac{n}{2}}}\int_{M}e^{-f}dv_{g}=1\},
$$
where $W(g,f,\tau)$ is Perelman's $W$ functional introduced in
\cite{perelman}. We recall that W-functional is defined by
$$
W(g,f,\tau)=(4\pi \tau)^{-n/2}\int_M[\tau(4|\nabla
u|^2+Ru^2)-2u^2\log u-nu^2)]dv_g
$$
with $ u=e^{-f/2}$, where
$$
\frac{1}{(4\pi\tau)^{\frac{n}{2}}}\int_{M}e^{-f}dv_{g}=1.
$$
It seems to us that a similar result should be true, but we can not
prove it yet.

On the spin manifold, we have the following result, which related
the first eigenvalue of Dirac operator to Perelman's lambda
constant.
\begin{thm}\label{thm4} Given a compact spin manifold with its Dirac operator
$D$. Let $\lambda^D_1$ be the first eigenvalue of $D$ in the sense
that the absolute value of $\lambda^D_1$ is the minimum of all
eigenvalues of $D$. Then we have
$${\lambda^D_1}^2\geq\lambda(g).$$
\end{thm}

Along the Ricci flow, L.Ma \cite{ma} obtained a interesting
monotonicity formula for eigenvalues of the Laplacian operators.

Here is the plan of the paper. in section two, we introduce some
standard formulae. We prove Theorem \ref{thm1} in section three.
Theorem \ref{thm2} is proved in section four and Theorem
\ref{thm3} is proved in section five. In the last section, we
discuss the relation between the first eigenvalue of the Dirac
operator and $\lambda$ constant of the $F(g,f)$ functional on spin
manifolds.
\section{preliminary}
In this section we will discuss the relationship between the
Laplace operator between the warped product space and its fibre.
The material below is standard (see \cite{Lo1}).
\begin{lem}\label {lem1}
Suppose that a Riemannian submersion $p:M\rightarrow B$ has
compact fibre F. For a smooth function $\phi^{M}$ on $M$, we let
$\phi^{B}=\int_{F}\phi^{M}dv^{F}$. Then we have the following
identity
$$
\Delta
^{B}\phi^{B}=\int_{F}\nabla^{2}\phi^{M}+[\breve{\delta}N+|\frac{\nabla_{hor}\phi^{M}}{\phi^{M}}
-N|^{2}-|\frac{\nabla_{hor}\phi^{M}}{\phi^{M}}|^{2}]\phi^{M}dv^{F}
$$
When the Riemannian submersion is a warped product
$M=B\times_{f}F,g^{M}=g^{B}+f^{2}g^{B}$, we have
$$
\Delta^{B}\phi^{B}=\int_{F}\nabla^{2}_{hor}\phi^{M}dv^{F}+p[\frac{\Delta^{B}f}{f}+\frac{|df|^2}{f^{2}}]\phi^{B}
-|N|^{2}\phi^{B}-2N\phi^{B}.
$$
\end{lem}
\begin{proof}
The first identity is from \cite{Lo2}. Assume now that $M$ is a
warped product. From \cite{besse}, we know that
$$
N=-\frac{p}{f}\nabla_{hor}f,
$$
and $$
\breve{\delta}N=p[\frac{\Delta^{B}f}{f}+\frac{|df|^{2}}{f^{2}}]
$$ where $p=dim F$.
Then,
\begin{eqnarray*}
\Delta^{B}\phi^{B}&=&\int_{F}\Delta^{B}\phi^{M}+[\breve{\delta}N
+|\frac{\nabla_{hor}\phi^{M}}{\phi^{M}}-N|^{2}
-|\frac{\nabla_{hor}}{\phi^{M}}|^{2}]\phi^{M}dv^{F}\\
&=&\int_{F}\Delta^{B}\phi^{M}dv^{F}+\int_{F}\breve{\delta}N\phi^{M}dv^{F}+\int_{F}|N|^{2}\phi^{M}dv^{F}\\
& &-2\int_{F}<\nabla_{hor}\phi^{M},N>dv^{F}
\end{eqnarray*}
Since in the warped product space, $\breve{\delta}N$, $|N|^{2}$
are constants on the fibre,we have
$$
\Delta^{B}\phi^{B}=\int_{F}\Delta^{B}\phi^{M}dv^{F}+(\breve{\delta}N+|N|^{2})\phi^{B}-2\int_{F}<\nabla_{hor}\phi^{M},N>dv^{F}
$$
Since we know (see \cite{Lo2}) that
$$
N\phi^{B}=\int_{F}(<\nabla_{hor}\phi^{M},N>-|N|^{2}\phi^{M})dv^{F},
$$
we have
\begin{eqnarray*}
\Delta^{B}\phi^{B}&=&\int_{F}\Delta^{B}\phi^{M}dv^{F}+(\breve{\delta}N-|N|^{2})\phi^{B}-2N\phi^{B}\\
&=&\int_{F}\Delta^{B}\phi^{M}dv^{F}+(p\frac{\Delta^{B}f}{f}+p(1-p)\frac{|df|^{2}}{f^{2}})\phi^{B}-2N\phi^{B}.
\end{eqnarray*}
\end{proof}
\section{proof of theorem \ref{thm1}}
Let $\lambda_1^M$ and $\phi^M$ be the first eigenvalue and the
eigenfunction of the Laplacian on $M$ respectively, ie., $$
\Delta^M\phi^M=-\lambda^M_1\phi^M.
$$
Similarly, we define $\lambda_1^B$, $\phi^B$, $\lambda_1^F$, and
$\phi^F$ respectively. Then
 for $u=\phi^B\phi^F$, we have
 $$
\Delta^Mu=-\lambda_1^Bu-\lambda_1^Fu-<\nabla_{hor}u,N>.
 $$
 Recall that
 $$
\bar{\lambda}_F=sup_{b\in B}\lambda_1^F(b).
 $$
By definition, we have
$$
\lambda_1^M=inf_{\{u\not=0\}}\frac{\int_M<-\Delta^Mu,u>}{\int_Mu^2}.
$$
Hence
$$
\lambda_1^M\leq \lambda_1^B+\bar{\lambda}_F+
\frac{\int_M<\nabla_{hor}u,N>u}{\int_Mu^2}.
$$
Define
$$
A_1=\frac{1}{2}sup_M(-\sum\nabla_{e_{\alpha}}<N,e_{\alpha}>).
$$
Then we have
$$
\lambda_1^M\leq \lambda_1^B+\bar{\lambda}_F+A_1.
$$
This is the upper bound for $\lambda_1^M$.

We now give a lower bound for $\lambda_1^M$. Note that
$$
-\Delta^M \phi^M=\lambda_1^M\phi^M.
$$
For simplicity, we let $u=\phi^M$. Then we have
$$
-\lambda_1^Mu=\nabla^2_{hor}u+\Delta^Fu-<\nabla_{hor}u,N>.
$$
Integrating over the fiber $F$, we have
$$
-\lambda_1^Mu^B=\int_{F}\nabla^{2}_{hor}u-\int_{F}<\nabla_{hor}u,N>dv^{F}
$$
By Lemma \ref{lem1}, we have
$$
-\lambda_1^Mu^B=\Delta^{B}u^{B}-(p\frac{\Delta^{B}f}{f}+p\frac{|df|^{2}}{f^{2}})\phi^{B}+Nu^{B}
$$

 Multiplying both sides of the equation
above by $u^B$ and taking integration on $B$, we get
\begin{eqnarray*}
\lambda_1^M\int_B|u^B|^2dv^{B}
&=&\int_B(|\nabla^Bu^B|^2dv^{B}-\frac{1}{2}\int_{B}Nu_{B}^{2}dv^{B}
\\
&+&\int_{B}p[\frac{\Delta^{B}
f}{f}+\frac{|df|^{2}}{f^{2}}]u_{B}^{2}dv^{B}\\
&=&\int_B|\nabla^Bu^B|^2dv^{B}-\frac{1}{2}\int_{B}<\nabla_{hor}
u_{B}^{2},-\frac{p}{f}\nabla_{hor}
f>dv^{B}\\
& &+\int_{B}p[\frac{\Delta^{B}
f}{f}+\frac{|df|^{2}}{f^{2}}]u_{B}^{2}dv^{B}\\
&=&\int_B|\nabla^Bu^B|^2dv^{B}-\frac{p}{2}\int_{B}u_{B}^{2}(\frac{\Delta^{B}
f}{f}-\frac{|\nabla
f|^{2}}{f^{2}})\\
& &+\int_{B}p[\frac{\Delta^{B}
f}{f}+\frac{|df|^{2}}{f^{2}}]u_{B}^{2}dv^{B}\\
&=&\int_B|\nabla^Bu^B|^2dv^{B}+\int_{B}(\frac{p}{2}\frac{\Delta^{B}
f}{f}+\frac{3p}{2}\frac{|df|^{2}}{f^{2}})u_{B}^{2}dv^{B}.
\end{eqnarray*}
Let
$$
c=inf (\frac{p}{2}\frac{\Delta^{B}
f}{f}+\frac{3p}{2}\frac{|df|^{2}}{f^{2}})
$$
Then we have
$$
\lambda_{1}^{M}\geq\lambda_{1}^{B}+c.
$$

\section{proof of theorem \ref{thm2}}
We shall use the same spirit of previous section to get bounds for
lambda constants.
\begin{proof}
 Suppose
$\phi^{M}$ is the first eigenfunction of the F-functional
$F(g,f)$. Then
\begin{eqnarray*}
\lambda^{M}_{1}\phi^{M}&=&-4\Delta ^{M}\phi^{M}+R^{M}\phi^{M}.
\end{eqnarray*}
Using the the same procedure in last section, we first do
integration in the fibre $F$.
\begin{eqnarray*}
\lambda_{1}^{M}\phi^{B}&=&-4\int_{F}\Delta
^{M}\phi^{M}dv^{F}+\int_{F}R^{M}\phi^{M}dv^{F}\\
\end{eqnarray*}
where $\phi^{B}=\int_{F}\phi^{F}dv^{F}$. From \cite{besse}, we
know
$$
R^{M}=R^{F}+R^{B}-p(p-1)\frac{|df|^{2}}{f^{2}}.
$$
Since the eigenfunction is positive, we have
\begin{eqnarray*}
\int_{F}R^{M}\phi^{M}dv^{F}&=&\int_{F}(R^{F}+R^{B}-p(p-1)\frac{|df|^{2}}{f^{2}})\phi^{M}dv^{F}\\
&=&(R^{B}-p(p-1)\frac{|df|^{2}}{f^{2}})\phi^{B}+\int_{F}R^{F}\phi^{M}dv^{F}\\
&\geq&(R^{B}-p(p-1)\frac{|df|^{2}}{f^{2}})\phi^{B}+\frac{1}{f^{2}}R_{0}\phi^{B}.
\end{eqnarray*}
Then, we multiply both sides by $\phi^{B}$, and integrate the
result on $B$. Using the result of the last theorem, we have
\begin{eqnarray*}
\lambda_{1}^{M}\int_{B}(\phi^{B})^{2}&=&-4\int_{B}\phi^{B}(\int_{F}\Delta^{M}\phi^{M}dv^{F})dv^{B}
+\int_{B}\phi^{B}\int_{F}R^{M}\phi^{M}dv^{F}dv^{B}\\
&=&4\int_B|\nabla^B\phi^B|^2dv^{B}+4\int_{B}(\frac{p}{2}\frac{\Delta^{B}f}{f}
+\frac{3p}{2}\frac{|df|^{2}}{f^{2}})\phi_{B}^{2}dv^{B}\\
&
&+\int_{B}\phi^{B}\int_{F}R^{B}+R^{F}-p(p-1)\frac{|df|^{2}}{f^{2}})\phi^{M}dv^{F}dv^{B}.\\
&=&\int_{B}4|\nabla^{B}\phi^{B}|^{2}+R^{B}(\phi^{B})^{2}dv^{B}\\
& &+\int_{B}(2p\frac{\Delta^{B}{f}}{f}+p(7-p)\frac{|df|^{2}}{f^{2}})(\phi^{B})^{2}dv^{B}\\
&
&+\int_{B}\phi^{B}\int_{F}\frac{1}{f^{2}}R_{0}\phi^{M}dv^{F}dv^{B}\\
&
&\geq\lambda_{1}^{B}\int_{B}(\phi^{B})^{2}dv^{B}+(c+\frac{1}{f^{2}}r_{0})\int_{B}(\phi^{B})^{2}dv^{B},
\end{eqnarray*}
where we have used the notation $$ c=inf
(2p\frac{\Delta^{B}f}{f}+p(7-p)\frac{|df|^{2}}{f^{2}}).$$
\end{proof}

\section{Related constants}
In this section, we shall discuss the relation between the first
eigenvalue of the Dirac functional and the lambda constants on
spin manifolds.

We now prove Theorem \ref{thm4}.

 \begin{proof} We consider the
Dirac functional
$$
I(u)=\int_M|Du|^2
$$
on the spin manifold $(M,g)$, where $u$ is the spinor field on $M$.
Using the Lichnorowicz formula
$$
D^2=\nabla^2+\frac{R}{4},
$$
we have
$$
I(u)=\int_M(|\nabla u|^2+\frac{R}{4}|u|^2)
$$
which is very similar to Perelman's F-functional. In fact, using
Kato's inequality
$$
|\nabla u|\geq |d|u||,
$$
we have
$$
I(u)\geq J(g,f)\geq \lambda(g)
$$
for $e^{-f}=|u|^2$. Let $\lambda^D_1$ be the first eigenvalue of
Dirac operator $D$, i.e.,
$$
Du_D=\lambda^D_1u_D.
$$
Then we have
$$
{\lambda^D_1}^2=I(u_D)\geq \lambda (g).
$$
\end{proof}
By this inequality, we can get a nice lower bound for
${\lambda^D_1}^2$ by studying $\lambda(g)$. We believe that this
kind of idea should be useful in the study of Dirac operators.

\end{document}